\documentclass[12pt]{article}

\usepackage{amsmath, epsfig, cite}
\usepackage{amssymb}
\usepackage{amsfonts}
\usepackage{latexsym}
\usepackage{amsthm}

\newtheorem{thm}{Theorem}[section]

\newtheorem{lem}[thm]{Lemma}

\numberwithin{equation}{section}

\renewcommand{\thefootnote}{}

\begin{document}

\begin{center}
{\large\bf On the  Askey--Wilson type integrals

 \footnote{ Email Address: weichuanan78@163.com. The work is supported by the National Natural Science Foundation of China (No. 12071103).}}
\end{center}

\renewcommand{\thefootnote}{$\dagger$}

\vskip 2mm \centerline{Chuanan Wei}
\begin{center}
{School of Biomedical Information and Engineering\\
  Hainan Medical University, Haikou 571199, China}
\end{center}


\vskip 0.7cm \noindent{\bf Abstract.} The Askey--Wilson integral is
very important in the theory of orthogonal polynomials. Liu's
integral is a generalization of the Askey--Wilson integral with many
parameters. With the help of the series rearrangement method, we
give the elementary proof of them. Furthermore, we establish two new
Askey--Wilson type integrals in the similar way and find a
generalization of a known transformation formula containing three
$_{3}\phi_{2}$ series.

\vskip 3mm \noindent {\it Keywords}: the Askey--Wilson integral;
   Liu's integral; the $q$-Gauss sum; the
nonterminating form of $q$-Vandermonde sum

 \vskip 0.2cm \noindent{\it AMS
Subject Classifications:} 05A30; 33D15.

\section{Introduction}
For two complex numbers $x$, $q$ with $|q|<1$, define the
$q$-shifted factorial to be
\begin{align*}
(x;q)_{\infty}=\prod_{k=0}^{\infty}(1-xq^k),\quad
(x;q)_n=\frac{(x;q)_{\infty}}{(xq^n;q)_{\infty}}\quad\text{with}\quad
 n\in\mathbb{Z}^{+}\cup\{0\}.
 \end{align*}
For convenience, we shall also adopt the following notation:
\begin{align*}
(x_1,x_2,\ldots,x_r;q)_{m}=(x_1;q)_{m}(x_2;q)_{m}\cdots(x_r;q)_{m},
 \end{align*}
 where $m\in\mathbb{Z}^{+}\cup\{0,\infty\}$.
Following Gasper and Rahman \cite{gasper}, define the basic
hypergeometric series by
$$
_{r+1}\phi_{r}\left[\begin{array}{c}
a_1,a_2,\ldots,a_{r+1}\\
b_1,b_2,\ldots,b_{r}
\end{array};q,\, z
\right] =\sum_{k=0}^{\infty}\frac{(a_1,a_2,\ldots, a_{r+1};q)_k}
{(q,b_1,b_2,\ldots,b_{r};q)_k}z^k.
$$
Then the $q$-binomial theorem (cf. \cite[Appendix (II. 3)]{gasper}),
the $q$-Gauss sum (cf. \cite[Appendix (II. 8)]{gasper}), and the
nonterminating form of $q$-Vandermonde sum (cf. \cite[Appendix (II.
23)]{gasper}) can be stated as

 \begin{align}\label{q-binomial}
 _{1}\phi_{0}\left[\begin{array}{c}
a\\-
\end{array};q,\, z
\right]=\frac{(az;q)_{\infty}}{(z;q)_{\infty}}\quad\text{with}\quad
 |z|<1,
 \end{align}
 \begin{align}
&_{2}\phi_{1}\left[\begin{array}{c} a, b\\c
\end{array};q,\, \frac{c}{ab}
\right]=\frac{(c/a,c/b;q)_{\infty}}{(c,c/ab;q)_{\infty}}
 \quad\text{with}\quad
 |c/ab|<1,
\label{q-gauss}
 \end{align}
 \begin{align}
&_{2}\phi_{1}\left[\begin{array}{c} a, b\\c
\end{array};q,\, q
\right] +\frac{(a,b,q/c;q)_{\infty}}{(qa/c,qb/c,c/q;q)_{\infty}}
{_{2}\phi_{1}}\left[\begin{array}{c} qa/c,qb/c\\q^2/c
\end{array};q,\, q
\right]
\notag\\[5pt]
&\quad=\frac{(qab/c,q/c;q)_{\infty}}{(qa/c,qb/c;q)_{\infty}}.
\label{basic}
  \end{align}

For $x=\cos\theta$, define the functions $h(x;\lambda)$ and
$h(x;\lambda_1,\lambda_2,\ldots,\lambda_s)$
 to be
\begin{align*}
 &h(x;\lambda)=(\lambda e^{i\theta},\lambda e^{-i\theta};q)_{\infty}
=\prod_{k=0}^{\infty}(1-2\lambda
xq^k+\lambda^2q^{2k}),\\[1mm]
&\:\:h(x;\lambda_1,\lambda_2,\ldots,\lambda_s)=h(x;\lambda_1)h(x;\lambda_2)\cdots
h(x;\lambda_s).
\end{align*}
 In 1985, Askey and Wilson \cite{askey} discovered an important
extension of the beta integral (see also \cite[Chapter 6]{gasper}),
which is now named the Askey--Wilson integral.

\begin{thm}\label{thm-a}
Let $a,b,c,d$ be complex numbers such that
$\max\{|a|,|b|,|c|,|d|\}<1$. Then
 \begin{align}\label{Askey}
\int_{0}^{\pi}\frac{h(x;1,-1,q^{\frac{1}{2}},-q^{\frac{1}{2}})}{h(x;a,b,c,d)}\,d\theta
 =\frac{2\pi(abcd;q)_{\infty}}{(q,ab,ac,ad,bc,bd,cd;q)_{\infty}}.
 \end{align}
\end{thm}

The original proof of Theorem \ref{thm-a} comes from the contour
integral method. More proofs of it can be seen in the papers
\cite{bowman,ismail-b,liu-a,liu-b,liu-d,rahman-a}.

 Recently,  Liu \cite{liu-c} established the following
parametric generalization of Theorem \ref{thm-a}.

\begin{thm}\label{thm-b}
Let $a,b,c,d,f,r,s,t,z,\alpha,\beta,\delta$ be complex numbers. Then
 \begin{align}\label{Liu}
&\int_{0}^{\pi}\frac{h(x;1,-1,q^{\frac{1}{2}},-q^{\frac{1}{2}})}{h(x;a,b,c,d,f)}
{_{4}\phi_{3}}\left[\begin{array}{c}ae^{i\theta},ae^{-i\theta},\beta,\delta\\r,s,t
\end{array};q,\, bcdfz \right]
d\theta
 \notag\\[1mm]
 &\:=\frac{2\pi(abcd,abcf,abdf,acdf;q)_{\infty}}{(q,ab,ac,ad,af,bc,bd,bf,cd,cf,df,q\alpha;q)_{\infty}}
  \notag\\[1mm]
 &\quad\times\:
  \sum_{n=0}^{\infty}\frac{1-\alpha q^{2n}}{1-\alpha}\frac{(\alpha,ab,ac,ad,af;q)_n}{(q,abcd,abcf,abdf,acdf;q)_n}
q^{\binom{n}{2}}(-bcdf)^n
 \notag\\[1mm]
 &\qquad\times
{_{4}\phi_{3}}\left[\begin{array}{c}q^{-n},\alpha
q^n,\beta,\delta\\r,s,t
\end{array};q,\, qz \right],
 \end{align}
where $\alpha=a^2bcdf/q$ and
$\max\{|a|,|b|,|c|,|d|,|f|,|r|,|s|,|t|,|z|\}<1$.
\end{thm}

 Theorem \ref{thm-b} includes some nice results as
special cases. Setting $r=a\mu,s=\beta,t=\delta,z=\mu/abcdf$ in
Theorem \ref{thm-b} and using \eqref{q-gauss}, Liu obtained
 \begin{align*}
\int_{0}^{\pi}\frac{h(x;1,-1,q^{\frac{1}{2}},-q^{\frac{1}{2}},\mu)}{h(x;a,b,c,d,f)}\,d\theta
 &=\frac{2\pi(\mu/a,a\mu,abcd,abcf,abdf,acdf;q)_{\infty}}{(q,ab,ac,ad,af,bc,bd,bf,cd,cf,df,a^2bcdf;q)_{\infty}}
 \notag\\[1mm]
 &\quad\times{_8W_7}(a^2bcdf/q;ab,ac,ad,af,abcdf/\mu;q,\mu/a),
 \end{align*}
where the symbol on the right-hand side denotes
   \begin{align*}
 \qquad _{8}W_{7}(c_1;c_2,\ldots,c_6;q,z)={_{8}\phi_{7}}\left[\begin{array}{c}
 c_1,q\sqrt{c_1},-q\sqrt{c_1},c_2,\ldots,c_6\\\sqrt{c_1},-\sqrt{c_1},qc_1/c_2,\ldots,qc_1/c_6
\end{array};q,\, z \right].
 \end{align*}
Under the transformation formula (cf. \cite[{Appendix (III.
23)}]{gasper}):
  \begin{align*}
  {_8W_7}(a;b,c,d,e,f;q,q^2a^2/bcdef)
  &=\frac{(qa,qa/ef,q\lambda/e,q\lambda/f;q)_{\infty}}{(qa/e,qa/f,q\lambda,q\lambda/ef;q)_{\infty}} \\[1mm]
  &\quad\times{_8W_7}(\lambda;\lambda b/a,\lambda c/a,\lambda d/a,e,f;q,qa/ef),
  \end{align*}
provided $\lambda=qa^2/bcd$, it is not difficult to understand that
the last conclusion is equivalent to the Nassrallah--Rahman integral
(cf. \cite[{Equation (6.3.8)}]{gasper}):
 \begin{align}
\int_{0}^{\pi}\frac{h(x;1,-1,q^{\frac{1}{2}},-q^{\frac{1}{2}},\mu)}{h(x;a,b,c,d,f)}\,d\theta
 &=\frac{2\pi(a\mu,b\mu,c\mu,abcd,abcf;q)_{\infty}}{(q,ab,ac,ad,af,bc,bd,bf,cd,cf,abc\mu;q)_{\infty}}
 \notag\\[1mm]
 &\quad\times{_8W_7}(abc\mu/q; ab,ac,bc,\mu/d,\mu/f;q,df).
\label{Nassrallah}
 \end{align}

In the same paper, he also displayed other two special cases of
Theorem \ref{thm-b}:
 \begin{align*}
&\int_{0}^{\pi}\frac{h(x;1,-1,q^{\frac{1}{2}},-q^{\frac{1}{2}})}{h(x;a,b,c,d,f)}
{_{3}\phi_{2}}\left[\begin{array}{c}ae^{i\theta},ae^{-i\theta},\alpha
rs/q\\\alpha r,\alpha s \end{array};q,\, bcdf \right] d\theta
 \\[1mm]
 &\:=\frac{2\pi(abcd,abcf,abdf,acdf;q)_{\infty}}{(q,ab,ac,ad,af,bc,bd,bf,cd,cf,df,a^{2}bcdf;q)_{\infty}}
  \\[1mm]
 &\quad\times
  \sum_{n=0}^{\infty}\frac{1-\alpha q^{2n}}{1-\alpha}\frac{(\alpha,ab,ac,ad,af,q/r,q/s;q)_n}{(q,abcd,abcf,abdf,acdf,\alpha r,\alpha
  s;q)_n}
q^{\binom{n}{2}}\bigg(-\frac{\alpha bcdfrs}{q}\bigg)^n,
 \\[2mm]
&\int_{0}^{\pi}\frac{h(x;1,-1,q^{\frac{1}{2}},-q^{\frac{1}{2}})}{h(x;a,b,c,d,f)}
{_{4}\phi_{3}}\left[\begin{array}{c}ae^{i\theta},ae^{-i\theta},\sqrt{r},-\sqrt{r}\\\sqrt{q\alpha},-\sqrt{q\alpha},r
\end{array};q,\, bcdf \right]
d\theta
 \\[1mm]
 &=\:\frac{2\pi(abcd,abcf,abdf,acdf;q)_{\infty}}{(q,ab,ac,ad,af,bc,bd,bf,cd,cf,df,a^{2}bcdf;q)_{\infty}}
 \\[1mm]
 &\quad\times\:
  \sum_{n=0}^{\infty}\frac{1-\alpha q^{4n}}{1-\alpha}\frac{(\alpha,ab,ac,ad,af;q)_{2n}(q,q\alpha/r;q^{2})_n}
  {(q,abcd,abcf,abdf,acdf;q)_{2n}(q\alpha,qr;q^{2})_n}
r^n(bcdf)^{2n}q^{2n^2-2n}.
 \end{align*}

 The rest of the paper is arranged as follows. We shall prove
Theorems \ref{thm-a} and \ref{thm-b} through the series
rearrangement method in Sections 2 and 3, respectively. Two new
Askey--Wilson type integrals and a symmetric formula involving
double series, which is
 a generalization of \eqref{trans-b}, will be derived in Section 4.

\section{An elementary proof of Theorem \ref{thm-a}}
In order to prove Theorem \ref{thm-a}, we need the following two
lemmas.

\begin{lem}\label{lemma-b}
 Let $c_1,\ldots,c_r,\lambda,\eta$ be
complex numbers. Then
 \begin{align*}
  \int_0^{\pi}\frac{h(x;\lambda)}{h(x;\eta)}\Psi(x;c_1,\ldots,c_r)d\theta
  =(\lambda \eta,\lambda/\eta;q)_{\infty}\sum_{k=0}^{\infty}\frac{(\lambda/\eta)^k}{(q,\lambda \eta;q)_k}
\int_0^{\pi}\frac{\Psi(x;c_1,\ldots,c_r)}{h(x;\eta q^k)}d\theta,
   \end{align*}
where $\Psi(x;c_1,\ldots,c_r)$ is an arbitrary expression associated
with $x$.
\end{lem}

\begin{proof}
By means of \eqref{q-gauss},
 there holds the following relation:
 \begin{align*}
 &\int_0^{\pi}\frac{\Psi(x;c_1,\ldots,c_r)}{h(x;\eta)}
  \frac{(\lambda e^{i\theta},\lambda e^{-i\theta};q)_{\infty}}{(\lambda\eta,\lambda/\eta;q)_{\infty}}d\theta
  \\[2mm]
  &=\int_0^{\pi}\frac{\Psi(x;c_1,\ldots,c_r)}{h(x;\eta)}
{_{2}\phi_{1}}\left[\begin{array}{c}\eta e^{i\theta},\eta
e^{-i\theta}\\\lambda\eta
\end{array};q,\,\frac{\lambda}{\eta} \right]
  d\theta
  \\[2mm]
  &=\sum_{k=0}^{\infty}\frac{(\lambda/\eta)^k}{(q,\lambda\eta;q)_k}
\int_0^{\pi}\!\frac{\Psi(x;c_1,\ldots,c_r)}{h(x;\eta q^k)}d\theta.
   \end{align*}
This finishes the proof of Lemma \ref{lemma-b}.
\end{proof}

\begin{lem}\label{lemma-a}
 Let $c_1,\ldots,c_r,\lambda,\eta$ be
complex numbers. Then
  \begin{align*}
  \int_0^{\pi}\frac{\Omega(x;c_1,\ldots,c_r)}{h(x;\lambda,\eta)}d\theta
  &=\frac{1}{(\lambda\eta,\eta/\lambda;q)_{\infty}}\sum_{k=0}^{\infty}\frac{q^k}{(q,q\lambda/\eta;q)_k}
\int_0^{\pi}\frac{\Omega(x;c_1,\ldots,c_r)}{h(x;\lambda q^k)}d\theta
\\[2mm]
&\quad+\text{idem}(\lambda;\eta),
   \end{align*}
where $\Omega(x;c_1,\ldots,c_r)$ is an arbitrary expression related
to $x$ and the notation $idem(\lambda;\eta)$ after an expression
signifies that the front expression is repeated with $\lambda$ and
$\eta$ interchanged.
\end{lem}

\begin{proof}
In terms of \eqref{basic}, we have
 \begin{align*}
  &\int_0^{\pi}\!\frac{\Omega(x;c_1,\ldots,c_r)}{h(x;\lambda)}
  \frac{(\lambda\eta,\eta/\lambda;q)_{\infty}}{(\eta e^{i\theta},\eta
  e^{-i\theta};q)_{\infty}}d\theta
  \\[2mm]
  &=\int_0^{\pi}\!\frac{\Omega(x;c_1,\ldots,c_r)}{h(x;\lambda)}
  {_{2}\phi_{1}}\left[\begin{array}{c}\lambda e^{i\theta},\lambda
e^{-i\theta}\\q\lambda/\eta
\end{array};q,\,q \right]d\theta
  \\[2mm]
&\quad+\int_0^{\pi}\!\frac{\Omega(x;c_1,\ldots,c_r)}{h(x;\lambda)}
  \frac{(\eta/\lambda,\lambda e^{i\theta},\lambda e^{-i\theta};q)_{\infty}}{(\lambda/\eta,\eta e^{i\theta},\eta
  e^{-i\theta};q)_{\infty}}
{_{2}\phi_{1}}\left[\begin{array}{c}\eta e^{i\theta},\eta
e^{-i\theta}\\\lambda\eta
\end{array};q,\,q \right]
  d\theta
  \\[2mm]
  &=\,\sum_{k=0}^{\infty}\frac{q^k}{(q,q\lambda/\eta;q)_k}
\int_0^{\pi}\!\frac{\Omega(x;c_1,\ldots,c_r)}{h(x;\lambda
q^k)}d\theta
\\[2mm]
&\quad+\:\frac{(\eta/\lambda;q)_{\infty}}{(\lambda/\eta;q)_{\infty}}
\sum_{k=0}^{\infty}\frac{q^k}{(q,q\eta/\lambda;q)_k}
\int_0^{\pi}\!\frac{\Omega(x;c_1,\ldots,c_r)}{h(x;\eta q^k)}d\theta.
   \end{align*}
This completes the proof of Lemma \ref{lemma-a}.
\end{proof}

Now we begin to prove Theorem \ref{thm-a}.

\begin{proof}
Via \eqref{q-binomial} and \eqref{q-gauss},
 it is ordinary to show that
 \begin{align*}
\int_{0}^{\pi}\frac{h(x;q^{\frac{1}{2}})}{h(x;a)}\,d\theta
 &=\frac{1}{2}\int_{-\pi}^{\pi}\frac{h(x;q^{\frac{1}{2}})}{h(x;a)}\,d\theta
 \\[1mm]
&=\frac{1}{2}\int_{-\pi}^{\pi}\sum_{k=0}^{\infty}\sum_{l=0}^{\infty}\frac{(q^{\frac{1}{2}}/a;q)_k(q^{\frac{1}{2}}/a;q)_l}{(q;q)_k(q;q)_l}a^{k+l}e^{i(k-l)\theta}d\theta
 \\[1mm]
 &=\frac{1}{2}\sum_{k=0}^{\infty}\sum_{l=0}^{\infty}\frac{(q^{\frac{1}{2}}/a;q)_k(q^{\frac{1}{2}}/a;q)_l}{(q;q)_k(q;q)_l}a^{k+l}\int_{-\pi}^{\pi}e^{i(k-l)\theta}d\theta
  \\[1mm]
 &=\pi\sum_{k=0}^{\infty}\frac{(q^{\frac{1}{2}}/a,q^{\frac{1}{2}}/a;q)_k}{(q,q;q)_k}a^{2k}
 \\[1mm]
 &=\pi\frac{(q^{\frac{1}{2}}a,q^{\frac{1}{2}}a;q)_{\infty}}{(q,a^{2};q)_{\infty}}.
  \end{align*}

Apply Lemma \ref{lemma-b} to it and then calculate the series on the
right-hand side by \eqref{q-gauss} to get
 \begin{align*}
\int_{0}^{\pi}\frac{h(x;q^{\frac{1}{2}},-q^{\frac{1}{2}})}{h(x;a)}\,d\theta
=\pi\frac{(q^{\frac{1}{2}},-q^{\frac{1}{2}};q)_{\infty}}{(q,a,-a;q)_{\infty}}.
  \end{align*}
Employ Lemma \ref{lemma-b} to this equation and then evaluate the
series on the right-hand side by \eqref{q-binomial} to derive
 \begin{align*}
\int_{0}^{\pi}\frac{h(x;1,q^{\frac{1}{2}},-q^{\frac{1}{2}})}{h(x;a)}\,d\theta
=\frac{2\pi}{(q,-a;q)_{\infty}}.
 \end{align*}
Exploit Lemma \ref{lemma-b} to the last equation and then compute
the series on the right-hand side by \eqref{q-binomial} to arrive at
 \begin{align*}
\int_{0}^{\pi}\frac{h(x;1,-1,q^{\frac{1}{2}},-q^{\frac{1}{2}})}{h(x;a)}\,d\theta
=\frac{2\pi}{(q;q)_{\infty}}.
 \end{align*}
Applying Lemma \ref{lemma-a} to it and then calculating the series
on the right-hand side by \eqref{basic}, we have
 \begin{align*}
\int_{0}^{\pi}\frac{h(x;1,-1,q^{\frac{1}{2}},-q^{\frac{1}{2}})}{h(x;a,b)}\,d\theta
=\frac{2\pi}{(q,ab;q)_{\infty}}.
  \end{align*}
Employing Lemma \ref{lemma-a} to this equation and then evaluating
the series on the right-hand side by \eqref{basic}, we deduce
\begin{align*}
\int_{0}^{\pi}\frac{h(x;1,-1,q^{\frac{1}{2}},-q^{\frac{1}{2}})}{h(x;a,b,c)}\,d\theta
=\frac{2\pi}{(q,ab,ac,bc;q)_{\infty}}.
 \end{align*}
Exploiting Lemma \ref{lemma-a} to the last equation and then
computing the series on the right-hand side by \eqref{basic}, we are
led to \eqref{Askey}.
\end{proof}

\section{An elemenary proof of Theorem \ref{thm-b}}

For proving Theorem \ref{thm-b}, we require the following two
transformation formulas (cf. \cite[{Exercises 2.22 and
3.6}]{gasper}):
\begin{align}
&\sum_{n=0}^{\infty}\frac{1-aq^{2n}}{1-a}\frac{(a,b,c,d,e;q)_n}{(q,qa/b,qa/c,qa/d,qa/e;q)_n}
\bigg(-\frac{q^2a^2}{bcde}\bigg)^nq^{\binom{n}{2}}
\notag\\[1mm]
&\:=\frac{(qa,qa/de;q)_{\infty}}{(qa/d,qa/e;q)_{\infty}}
{_{3}\phi_{2}}\left[\begin{array}{c}qa/bc,d,e\\qa/b,qa/c
\end{array};q,\,\frac{qa}{de} \right],
\label{trans-a}
\end{align}
where $|aq/de|<1$,

\begin{align}
&{_{3}\phi_{2}}\left[\begin{array}{c}a,b,c\\d,e
\end{array};q,\,q \right]
+\frac{(a,b,c,qd/e,q/e;q)_{\infty}}{(qa/e,qb/e,qc/e,d,e/q;q)_{\infty}}
{_{3}\phi_{2}}\left[\begin{array}{c}qa/e,qb/e,qc/e\\qd/e,q^2/e
\end{array};q,\,q \right]
\notag\\[1mm]
&\:=\frac{(qab/e,qac/e,d/a,q/e;q)_{\infty}}{(qa/e,qb/e,qc/e,d;q)_{\infty}}
{_{3}\phi_{2}}\left[\begin{array}{c}a,qa/e,qabc/de\\qab/e,qac/e
\end{array};q,\,\frac{d}{a} \right],
\label{trans-b}
\end{align}
where $|d/a|<1$.

Now we start to prove Theorem \ref{thm-b}.

\begin{proof}
Let $P(a,b,c,d,f,r,s,t,z,\alpha,\beta,\delta)$ denotes the double
sum on the right-hand side of \eqref{Liu}. It is routine to verify
that
\begin{align*}
&P(a,b,c,d,f,r,s,t,z,\alpha,\beta,\delta)
\\[1mm]
 &=\:\sum_{n=0}^{\infty}\sum_{k=0}^{n}\frac{1-\alpha
q^{2n}}{1-\alpha}\frac{(\alpha,ab,ac,ad,af;q)_n}{(q,abcd,abcf,abdf,acdf;q)_n}
\frac{(q^{-n},\alpha
 q^n,\beta,\delta;q)_k}{(q,r,s,t;q)_k}
 \\[1mm]
 &\quad\:\times q^{\binom{n}{2}}(-bcdf)^n(qz)^k
 \\[1mm]
&=\:\sum_{k=0}^{\infty}\sum_{n=k}^{\infty}\frac{1-\alpha
q^{2n}}{1-\alpha}\frac{(\alpha,ab,ac,ad,af;q)_n}{(q,abcd,abcf,abdf,acdf;q)_n}
\frac{(q^{-n},\alpha
 q^n,\beta,\delta;q)_k}{(q,r,s,t;q)_k}
 \\[1mm]
 &\quad\:\times q^{\binom{n}{2}}(-bcdf)^n(qz)^k.
 \end{align*}

Shifting the index $n\to n+k$, we obtain
\begin{align*}
& P(a,b,c,d,f,r,s,t,z,\alpha,\beta,\delta)
\\[1mm]
&=\sum_{k=0}^{\infty}\frac{(ab,ac,ad,af,\beta,\delta;q)_k(q\alpha;q)_{2k}}{(q,abcd,abcf,abdf,acdf,r,s,t;q)_k}(bcdfz)^k
\\[1mm]
 &\:\quad\times\sum_{n=0}^{\infty}\frac{1-\alpha q^{2k+2n}}{1-\alpha q^{2k}}
 \frac{(\alpha q^{2k},abq^{k},acq^{k},adq^{k},afq^{k};q)_n}{(q,acdfq^{k},abdfq^{k},abcfq^{k},abcdq^{k};q)_n}q^{\binom{n}{2}}(-bcdf)^n
\end{align*}
\begin{align*}
 &=\frac{(q\alpha,df;q)_{\infty}}{(abdf,acdf;q)_{\infty}}\sum_{k=0}^{\infty}\frac{(ab,ac,ad,af,\beta,\delta;q)_k}{(q,abcd,abcf,r,s,t;q)_k}(bcdfz)^k
\\[1mm]
 &\:\quad\times
{_{3}\phi_{2}}\left[\begin{array}{c}abq^k,acq^k,bc\\abcdq^k,abcfq^k
\end{array};q,\,df \right].
\end{align*}
where we have utilized \eqref{trans-a} in the last step. Thus there
is the relation:
\begin{align}
&\sum_{k=0}^{\infty}\frac{(ab,ac,ad,af,\beta,\delta;q)_k}{(q,abcd,abcf,r,s,t;q)_k}(bcdfz)^k
{_{3}\phi_{2}}\left[\begin{array}{c}abq^k,acq^k,bc\\abcdq^k,abcfq^k
\end{array};q,\,df \right]
 \notag\\[1mm]
&=\frac{(abdf,acdf;q)_{\infty}}{(q\alpha,df;q)_{\infty}}\sum_{n=0}^{\infty}\frac{1-\alpha
q^{2n}}{1-\alpha}\frac{(\alpha,ab,ac,ad,af;q)_n}{(q,abcd,abcf,abdf,acdf;q)_n}
q^{\binom{n}{2}}(-bcdf)^n
 \notag\\[1mm]
 &\quad\times
{_{4}\phi_{3}}\left[\begin{array}{c}q^{-n},\alpha
 q^n,\beta,\delta\\r,s,t
\end{array};q,\,qz \right].
 \label{trans-c}
 \end{align}

 Applying Lemma \ref{lemma-a} to \eqref{Askey}, we have
 \begin{align}
\int_{0}^{\pi}\frac{h(x;1,-1,q^{\frac{1}{2}},-q^{\frac{1}{2}})}{h(x;a,b,c,d,f)}\,d\theta
 &=\frac{2\pi(abcd;q)_{\infty}}{(q,ab,ac,ad,bc,bd,cd,df,f/d;q)_{\infty}}
 \notag\\[1mm]
 &\quad\times
{_{3}\phi_{2}}\left[\begin{array}{c}ad,bd,cd\\qd/f,abcd
\end{array};q,\,q \right]
 +\text{idem}(d;f).
 \label{Integral}
\end{align}
Performing the replacements $a\to ad$, $b\to bd$,  $c\to cd$, $d\to
abcd$, $e\to qd/f$ in \eqref{trans-b}, there holds
\begin{align}
&\frac{(af,bf,cf,abcd;q)_{\infty}}{(bc,f/d,abdf,acdf;q)_{\infty}}{_{3}\phi_{2}}\left[\begin{array}{c}ad,bd,cd\\qd/f,abcd
\end{array};q,\,q \right]
+\text{idem}(d;f)
\notag\\[1mm]
&\:\:={_{3}\phi_{2}}\left[\begin{array}{c}ad,af,df\\abdf,acdf
\end{array};q,\,bc \right].
\label{trans-c}
\end{align}
Dealing with the expression on the right-hand side of
\eqref{Integral} by \eqref{trans-c} and then using the substitutions
$b\to d, c\to f, d\to c, f\to b$ in the resulting identity, we get
the Ismail--Stanton--Viennot integral (cf. \cite[Theorem
3.5]{ismail-c}):
 \begin{align}
\int_{0}^{\pi}\frac{h(x;1,-1,q^{\frac{1}{2}},-q^{\frac{1}{2}})}{h(x;a,b,c,d,f)}\,d\theta
 &=\frac{2\pi(abcd, abcf;q)_{\infty}}{(q,ab,ac,ad,af,bc,bd,bf,cd,cf;q)_{\infty}}
  \notag\\[1mm]
 &\times
{_{3}\phi_{2}}\left[\begin{array}{c}ab,ac,bc\\abcd,abcf
\end{array};q,\,df \right].
 \label{ismail}
\end{align}
Thanks to \eqref{trans-c} and \eqref{ismail}, we catch hold of
\eqref{Liu}.
\end{proof}

\section{Two new Askey-Wilson type integrals}

Above all, we shall display the following Askey-Wilson type
integral.

\begin{thm}\label{thm-d}
Let $a,b,c,d,f,g,\mu$ be complex numbers. Then
 \begin{align}
&\int_{0}^{\pi}\frac{h(x;1,-1,q^{\frac{1}{2}},-q^{\frac{1}{2}},\mu)}{h(x;a,b,c,d,f,g)}\,d\theta
 =\frac{2\pi(c\mu,d\mu,f\mu,acdf,bcdf;q)_{\infty}}{(q,ac,ad,af,bc,bd,bf,cd,cf,df,fg,g/f,cdf\mu;q)_{\infty}}
\notag \\[1mm]
&\quad\times\:\sum_{n=0}^{\infty}\frac{(af,bf,cf,df,cdf\mu;q)_{n}\,q^{n}}{(q,qf/g,f\mu,acdf,bcdf;q)_n}
{_8W_7}(cdf\mu q^{n-1};\mu/a,\mu/b,cd,cfq^{n},dfq^{n};q,ab)
\notag \\[1mm]
  &\quad+\:\,\text{idem}(f;g),
  \label{eq:wei-a}
 \end{align}
where $\max\{|a|,|b|,|c|,|d|,|f|,|g|\}<1$.
\end{thm}

\begin{proof}
Perform the replacements $a\to d, b\to f, d\to a,f\to b$ in
\eqref{Nassrallah} to achieve
\begin{align*}
\int_{0}^{\pi}\frac{h(x;1,-1,q^{\frac{1}{2}},-q^{\frac{1}{2}},\mu)}{h(x;a,b,c,d,f)}\,d\theta
 &=\frac{2\pi(c\mu,d\mu,f\mu,acdf,bcdf;q)_{\infty}}{(q,ac,ad,af,bc,bd,bf,cd,cf,df,cdf\mu;q)_{\infty}}
 \notag \\[1mm]
 &\quad\times\,{_8W_7}(cdf\mu/q;\mu/a,\mu/b,cd,cf,df;q,ab).
\end{align*}
Exploiting Lemma \ref{lemma-a} to the last equation, we obtain
\eqref{eq:wei-a} after some simplification.
\end{proof}

Subsequently, applying Lemma \ref{lemma-b} to \eqref{eq:wei-a}, we
get the following theorem after some simplification.

\begin{thm}\label{thm-e}
Let $a,b,c,d,f,g,\mu,\nu$ be complex numbers. Then
 \begin{align*}
&\int_{0}^{\pi}\!\frac{h(x;1,-1,q^{\frac{1}{2}},-q^{\frac{1}{2}},\mu,\nu)}{h(x;a,b,c,d,f,g)}\,d\theta
 =\frac{2\pi(c\mu,d\mu,d\nu,\nu/d,f\mu,acdf,bcdf;q)_{\infty}}{(q,ac,ad,af,bc,bd,bf,cd,cf,df,fg,g/f,cdf\mu;q)_{\infty}}
 \notag \\[1mm]
&\quad\times\:\sum_{m=0}^{\infty}\sum_{n=0}^{\infty}\frac{(ad,bd,cd;q)_{m}}{(q,d\mu,d\nu;q)_m}
\frac{(af,bf,cf;q)_{n}}{(q,qf/g,f\mu;q)_n}\frac{(df,cdf\mu;q)_{m+n}}{(acdf,bcdf;q)_{m+n}}\bigg(\frac{\nu}{d}\bigg)^{m}q^{n}
\notag \\[1mm]
&\:\qquad\times{_8W_7}(cdf\mu
q^{m+n-1};\mu/a,\mu/b,cdq^{m},cfq^{n},dfq^{m+n};q,ab)+\text{idem}(f;g),
 \end{align*}
where $\max\{|a|,|b|,|c|,|d|,|f|,|g|,|v/d|\}<1$.
\end{thm}

Finally, we shall provide the following symmetric formula.

\begin{thm}\label{thm-f}
Let $S(a,b,c,d;f,g)$ stand for the double series:
 \begin{align*}
&S(a,b,c,d;f,g)=\frac{(f/g,acdf,bcdf;q)_{\infty}}{(af,bf,df,c^2df;q)_{\infty}}
\notag \\[1mm]
&\quad\times\:\sum_{n=0}^{\infty}\frac{(af,bf,df,c^2df;q)_{n}\,q^{n}}{(q,qf/g,acdf,bcdf;q)_n}
{_8W_7}(c^2dfq^{n-1};c/a,c/b,cd,cfq^{n},dfq^{n};q,ab).
 \end{align*}
Then we have
\begin{align}
&S(a,b,c,d;f,g)+S(a,b,c,d;g,f)
\notag \\[2mm]
&\:=\frac{(ac,bc,bd,f/g,g/f,abfg,adfg;q)_{\infty}}
{(c^2,ab,af,ag,bf,bg,df,dg;q)_{\infty}}
{_{3}\phi_{2}}\left[\begin{array}{c}af,ag,fg\\abfg,adfg
\end{array};q,\,bd \right],
\label{symmetic}
 \end{align}
where $\max\{|ab|,|bd|\}<1$.
\end{thm}

\begin{proof}
Fix $\mu=c$ in \eqref{eq:wei-a} to deduce
\begin{align}
&\int_{0}^{\pi}\frac{h(x;1,-1,q^{\frac{1}{2}},-q^{\frac{1}{2}})}{h(x;a,b,d,f,g)}\,d\theta
 =\frac{2\pi(c^2,acdf,bcdf;q)_{\infty}}{(q,ac,ad,af,bc,bd,bf,df,fg,g/f,c^2df;q)_{\infty}}
\notag \\[1mm]
&\quad\times\:\sum_{n=0}^{\infty}\frac{(af,bf,df,c^2df;q)_{n}\,q^{n}}{(q,qf/g,acdf,bcdf;q)_n}
{_8W_7}(c^2dfq^{n-1};c/a,c/b,cd,cfq^{n},dfq^{n};q,ab)
\notag \\[1mm]
  &\quad+\:\,\text{idem}(f;g).
  \label{eq:wei-b}
 \end{align}
Utilizing the substitutions $b\to f$, $c\to g$,  $f\to b$ in
\eqref{ismail}, there is
\begin{align}
\int_{0}^{\pi}\frac{h(x;1,-1,q^{\frac{1}{2}},-q^{\frac{1}{2}})}{h(x;a,b,d,f,g)}\,d\theta
&= \frac{2\pi(abfg,adfg;q)_{\infty}}
{(q,ab,ad,af,ag,bf,bg,df,dg,fg;q)_{\infty}}
\notag \\[1mm]
&\times\,{_{3}\phi_{2}}\left[\begin{array}{c}af,ag,fg\\abfg,adfg
\end{array};q,\,bd \right].
  \label{eq:wei-c}
 \end{align}
Therefore, the combination of \eqref{eq:wei-b} and \eqref{eq:wei-c}
leads to \eqref{symmetic}.
\end{proof}

\textbf{Remark}

 Noting that \eqref{symmetic} reduces to
\eqref{trans-c} when $c=a$, so we can regard \eqref{symmetic} as a
generalization of \eqref{trans-b}.


\end{document}